\documentclass[twoside,a4paper,12pt]{article}
\usepackage{amsfonts, amsbsy, amsmath, amsthm, amssymb, latexsym, verbatim, enumerate}

\newtheorem{propo}{\sc Proposition}[section]
\newtheorem{prop}{\sc Proposition}[section]

\newtheorem{cor}[propo]{\sc Corollary}
\newtheorem{lem}[propo]{\sc Lemma}
\newtheorem{thm}[propo]{\sc Theorem}
\usepackage{hyperref}

\newcommand{\cupdot}{\mathbin{\mathaccent\cdot\cup}}

\newcommand{\pf}{{\it Proof:\quad}}

\newcommand{\ga}{\alpha}

\newcommand{\Ga}{\Gamma}

\newcommand{\gl}{\lambda}

\newcommand{\cp}{\mathbb{C}}

\newcommand{\dne}{\hfill $\Box$ \vspace{0.3cm}}

\begin{document}

\title{Singular  graphs }
\author{ Ali Sltan AL-Tarimshawy \
{\small Version 10 July 2017, printed \today}}
\maketitle

\begin{abstract}
Let $\Gamma$ be a simple  graph on a finite vertex set $V$ and let $A$ be its adjacency matrix. Then $\Gamma$ is said to be singular if and only if $0$ is an eigenvalue of $A.$ The nullity (singularity) of $\Gamma,$ denoted by ${\rm null}(\Gamma),$ is the algebraic multiplicity of the eigenvalue $0$ in the spectrum of $\Gamma.$ In 1957, Collatz and Sinogowitz \cite{von1957spektren} posed the problem of characterizing singular graphs. Singular graphs have important applications in mathematics and science. The chemical importance of singular graphs lies in the fact that if the nullity for the molecular graph is greater than zero then the corresponding chemical compound is highly reactive or unstable.  By this reason, the chemists have a great interest in this problem.  The general problem of characterising singular graphs is easy to state but it seems too difficult. In this work, we investigate this problem for  graphs in general and graphs with a vertex  transitive group $G$  of automorphisms. In some cases we determine the nullity of such graphs. We characterize singular Cayley graphs over cyclic  groups. We show that vertex transitive graphs with $|V|$ is prime are non-singular.  
\end{abstract}

\section{\sc Introduction}
Let $\Gamma$ be a  graph on the finite  vertex set $V$ of size $n,$ and let $A$ be its adjacency matrix. Then $\Gamma$ is \textit{singular} if ~$A$ is singular. The \textit{spectrum }of $\Gamma$  consists of all eigenvalues $\lambda_1,...,\lambda_n$ of $A$ and so $\Gamma$ is singular if and only if $0$ belongs to the spectrum of $\Gamma.$  The \textit{nullity} of  $\Gamma$ is the dimension of the null space of $\Gamma$ and  we denote this by ${\rm null}(\Gamma).$   Note $|V|={\rm null}(\Gamma) +r(\Gamma)$ where $r(\Gamma)$ is the rank  of $A.$   Hence  singular graphs have a non-trivial null space.

  There are applications of graph spectra and singularity in the representation theory of permutation groups. In physics and chemistry, nullity is important for the study of molecular graphs stability.  The nullity of a graph is also important in mathematics since it is relevant to the rank of the adjacency matrix.  Collatz and Sinogowitz \cite{von1957spektren} posed the problem of characterizing all graphs with zero nullity. These are the non-singular graphs.
  
  In this work, we determine conditions for a graph to be singular. In Section~4 we investigate this problem for Cayley graphs. In the last Section, we provide a new method for computing the spectrum of vertex transitive graphs and by this method we discuss the singularity of such graphs. 
  
  \section{\sc Representations and Characters of Groups} 
    
 Let $\mathbb{C}$ be the field of complex numbers.  Let $G$ be a finite group and $W$ be a finite dimensional  vector space over  $\mathbb{C}.$ A \textit{representation} of  $G$ over $\mathbb{C}$ is a group homomorphism $\rho$ from $G$ to $GL(W)$. Here $GL(W)$  is the group of all bijective linear maps $\beta:W\rightarrow W.$  The \textit{degree} of the representation $\rho$  is the dimension of the vector space $W.$  We also say that $\rho$ is a representation of degree $n$ over $\mathbb{C}.$ So if $\rho$ is a map from $G$ to $GL(W)$  then $\rho$ is a representation if and only if
 $$\rho(gh)=\rho(g)\rho(h) $$ for all $g,h \in G.$
 
  The \textit{character} associated  with $ \rho$  is the  function $ \chi_\rho:G\rightarrow \mathbb{C}$ denoted by $ \chi_\rho(g)={\rm tra}(\rho(g)) $ for all  $ g\in G$.  Here ${\rm tra}(\rho(g)) $ is the trace of the representation matrix. \textit{The degree} of the character is the degree of the representation and it is equal to $\chi_\rho(1_G)$. It is clear that characters are \textit{class functions} (functions which are constant on all conjugate classes) and the set of all irreducible characters is the  basis of  the vector space of all   class functions  on  $G.$

   Let $W=\mathbb{C}G$ be the vector space with basis $G$ see Section~3 for the definition. The \textit{right regular representation} is  a map (in fact homomorphism) $\rho_r:G\rightarrow GL(W)$ of $G$ given by $\rho_r(h)(g)=gh$ for each $h\in G$ and all $g\in G$. The \textit{left regular representation} $\rho_l:G\rightarrow GL(W)$ of $G$ is given by $\rho_l(h)(g)=h^{-1}g$ for each $h\in G$ and all $g\in G.$
   
    A $\mathbb{C}G$-\textit{module}  is a vector space $W$ over $\mathbb{C} $ if an action  $(w,g)\rightarrow w^g\in W(w\in W, g\in G)$ is defined satisfying the following conditions:
  
	(1) $w^g \in W$
	
	(2) $(w^g)^h=w^{gh}$
	
	(3) $(\lambda w)^g=\lambda w^g$
	
	 (4) $(u+w)^g=u^g+w^g$

for all $u,w\in W, \lambda\in \mathbb{C}$ and $g,h\in G.$ A subset $\overline{W}$ of $W$ is said to be an $\mathbb{C}G$-\textit{sub-module} of $W$ if $\overline{W}$ is a subspace and $\overline{w}^g\in \overline{W}$ for all $g\in G$ and for all $\overline{w}\in \overline{W}.$   An $\mathbb{C}G$-module $W$ is said to be \textit{irreducible} if it is  non-zero and it has no $\mathbb{C}G$-sub-modules other than $\{0\}$ and $W.$

Now we list some results that we need these  in this work:

\begin{thm}(Schur's Lemma in the terms of representations)\cite{james2001representations}\label{noor4}
Let $V$ and $W$ be irreducible $\mathbb{C}G$-modules.

(1) If $\varphi:V\longrightarrow W$ is a $\mathbb{C}G$-homomorphism then either $\varphi$ is a $\mathbb{C}G$-isomorphism or $\varphi(v)=0$ for all $ v\in V$.

(2) If $\varphi:V\longrightarrow V$ is a $\mathbb{C}G$-isomorphism then $\varphi=\lambda id_V$ where $\lambda\in \mathbb{C}$.
\end{thm}

 \begin{thm}(Orbit-Stabilizer Theorem) \label{zu6}\cite[p.~21]{godsil2013algebraic}
Let $G$ be a permutation group acting on a set $V$ and let $u$ be a point in $V$. Then $|G|=|G(u)||G_u|.$
\end{thm}

\section{\sc Graphs and associated Permutation Modules}

Let $\Ga=(V,E)$ be a graph with vertex set $V=\{v_{1},\,v_{2},\,...,\,v_{n}\}$ and edge set $E.$ We write also $V=V(\Ga)$ and $E=E(\Ga).$ All graphs in this paper are undirected without loops or multiple edges. Let $\cp V$ denote the vector space with basis $V$ over $\cp.$  Its elements are the formal sums $$f=\sum_{v\in V}\,f_{v}\,v$$ with $f_{v}\in \cp.$ In particular, we identify $v=1v$ so that $V$ is a subset of $\cp V.$ Equivalently, $\cp V$ may be regarded as the set of all functions $f\!:\, V\to \cp$ where $f$ is the function with $f(v)=f_{v}$ for all $v\in V.$ Occasionally we identify $f$ also with its coordinate vector in $\cp^{n}.$ This vector space has a natural inner product, given by  $\langle v,v'\rangle =1$ if $v=v'$ and $\langle v,v'\rangle=0$ if $v\neq v',$ for all $v,\,v'\in V.$ 

Two distinct vertices $v,\,v'$ are adjacent, denoted $v\sim v',$ if and only if $\{v,v'\}$ is an edge. The {\it adjacency map} $\alpha\!:\, \cp V\to \cp V$ is the linear map given by $$(*)\quad\alpha(v)=\sum_{v\sim v'}\,v'\,$$ for all $v\in V.$ 
Since $v\sim v'$ for $v,\,v'\in V$ if and only if $v'\sim v$ we have $\langle  \alpha(v),v'\rangle =\langle v,\alpha(v')\rangle.$ Therefore $\alpha$ is symmetric for the given  inner product.  The matrix of $\alpha$ with respect to the basis $V$ is  the {\it adjacency matrix} $A=A(\Ga)$ of $\Ga.$ 

Since $A$ is symmetric all eigenvalues of $A$ are real. We denote the distinct eigenvalues by  $\gl_{1}> \gl_{2}> ...> \gl_{t}$ and let $\mu_{1},\,\mu_{2},\,...,\mu_{t}$ be their multiplicities. We denote by $E_{1},\,E_{2},\,...\,E_{t}$ the corresponding eigenspaces. It follows from the symmetry of $A$ that  $\cp V$ decomposes into eigenspaces $$\cp V\,=\,E_{1}\,\oplus\, E_{2}\,\oplus\,...\,\oplus E_{t}$$
and that the multiplicity $\mu_{i}$ of $\gl_{i}$ is the dimension of $E_{i}.$ The
 {\it spectrum} of $\Ga$ consists of the eigenvalues $\gl_{1}^{\mu_{1}},\gl_{2}^{\mu_{2}},\,...,\gl_{t}^{\mu_{t}}$ of $A$ where $\gl_{1}^{\mu_{1}}$ indicates that $\gl_{1}$ has multiplicity $\mu_{1},$  and so on.  Throughout we denote the kernel of $\ga$ by $E_{*}.$ Thus $E_{*}=E_{\gl}$ when $\lambda=0$ is an eigenvalue of $\ga$ and  $E_{*}=0$ otherwise. Note  the {\it nullity} of  $\Ga$ is the dimension of $E_{*},$ it is denoted by ${\rm null}(\Ga)=\dim(E_{*}).$
   \begin{lem}\cite[p.~166]{godsil2013algebraic}\label{ma2}
Let $W$ and $U$ be  vector spaces with linear maps $$\varphi : W\rightarrow U \,\, \text{and}\,\, \varsigma :U \rightarrow W.$$ Then $\varphi \varsigma :U \rightarrow U$ and $\varsigma \varphi :W \rightarrow W$ have  the same nonzero eigenvalues. Furthermore, if $\lambda$ is a nonzero eigenvalue with eigenspace $W_\lambda \subseteq W$ and $U_\lambda\subseteq U$ for $\varsigma \varphi$ and $\varphi \varsigma$ respectively then $\varphi$ and $\varsigma$ restrict to isomorphism $\varphi : W_\lambda \rightarrow U_\lambda $ and $\varsigma:U_\lambda \rightarrow W_\lambda $. 
\end{lem}
Automorphisms of $\Ga$ are permutations $g$ of $V,$ denoted $g\!:\,v\to v^{g},$ such that $v\sim v'$ if and only if $v^{g}\sim v'^{g}$ for all $v,\,v'\in V.$ The group of all automorphisms of $\Ga$ is denoted by ${\rm Aut}(\Ga).$ The action of $g$ on $V$ can be extended to a linear map $g\!:\,\cp V\to \cp V$ by setting $$f^{g}=\sum_{v\in V}\,f_{v}\,v^{g}.$$ 
This turns $\cp V$ into a permutation module for ${\rm Aut}(\Ga)$ and we may  refer to $\cp V$ as the {\it vertex module} of $\Ga.$ The character of the permutation module is the number of fixed points of permutation. We  denote this  by  $\psi$, so $\psi(g)=$ is the number of $v_i$ such that $v_i^g=v_i$. The following is an essential property of the adjacency map. 

\medskip
\begin{prop}\label{zu}\cite[p.~134]{cvetkovic1995spectra}
A permutation $g$ of $V$ is an automorphism of $\Gamma$ if and only if $\alpha (f^g)= (\alpha(f))^g $ where $\alpha$ is the adjacency map of  $\Gamma$ and $f\in \mathbb{C}V.$ 
\end{prop}
\pf It suffices to show this property when $f=v$ for some $v\in V.$  Assume that $\alpha (v^g)= (\alpha(v))^g.$    Suppose that $u \sim v$ where $u,v\in V.$ Therefore we want to prove that $g$ is an automorphism of $\Gamma.$ So by the definition of the automorphism this  is enough to prove that $u^g \sim v^g.$ Since $u \sim v$ so $\langle\alpha(u),v\rangle =1.$ Then  
\begin{eqnarray}
\langle \alpha(u^g),v^g\rangle &=& \langle (\alpha(u))^g,v^g\rangle \nonumber\\
&=& \langle\alpha(u),v^{gg^{-1}}\rangle \nonumber\\
&=& \langle\alpha(u),v\rangle \nonumber\\
&=&1. \nonumber
\end{eqnarray}
Hence $u^g \sim v^g$ which means  that  $g$ is an automorphism of $\Gamma.$

Now suppose that $g$ is an automorphism of $\Gamma$ and we want to prove that $\alpha(v^g)=(\alpha(v))^g$ for all $v \in V.$ Then 
\begin{equation}\label{eigen12}
\alpha(v^g)=\sum_w\langle \alpha(v^g),w\rangle w
\end{equation}
\begin{equation}\label{eigen13}
(\alpha(v))^g=\sum_w\langle (\alpha(v))^g,w\rangle w.
\end{equation}
Hence $\langle \alpha(v^g),w\rangle =\langle v^g,\alpha(w)\rangle =   \left\lbrace
\begin{array}[pos]{cc}
1& v^g \sim w\\
0& v^g \nsim w\\
\end{array}
\right.
$

and $\langle  (\alpha(v))^g,w \rangle =\langle \alpha(v),w^{g^{-1}} \rangle =\left\lbrace
\begin{array}[pos]{cc}
1& v \sim w^{g^{-1}}\\
0& v \nsim w^{g^{-1}}.\\
\end{array}
\right.
$

Now, by the definition of  automorphisms, $v^g\sim w$ if and only if $v^{gg^{-1}}=v \sim w^{g^{-1}}$. Therefore by    \autoref{eigen12} and   \autoref{eigen13}  we have that $\alpha (f^g)=(\alpha(f))^g$ for all $ f \in \mathbb{C}V.$
\dne

\medskip

\begin{thm}\cite{biggs1993algebraic}\label{noor47}
Let $\Gamma$ be a finite graph with   adjacency map $\alpha$ and eigenspaces $E_1,E_2,\cdots,E_t$  corresponding to the distinct eigenvalues of $\alpha.$  Let $G$ be a group of automorphisms of $\Gamma.$  Then  each  $E_i$  is a $\mathbb{C}G$-module.
\end{thm}

\medskip   
The {\it complement\,} of the graph $\Ga=(V,E)$ is the graph $\overline{\Ga}$ on the same vertex set in which any two vertices $v,\, v'$ are adjacent to each other if and only if they are not adjacent to each other in $\Ga.$    The \textit{line graph} of  $\Gamma$ denoted by $L(\Gamma)$, is the graph whose vertex set is the edge set of $\Gamma.$  Two vertices are  adjacent in $L(\Gamma)$ if and only if these edges are incident in  $\Gamma$ (that is, the two edges have the same endpoint). The singularity of a graph is closely related to another special eigenvalue of graphs: 

\medskip
\begin{prop}\label{sing1} Let $\Ga$ be a regular graph. Then $\Ga$ is singular if and only if $-1$ is an eigenvalue of  $\overline{\Ga}.$
\end{prop}

\pf Let $|V|=n,$ let $I$ be the $n \times n$ identity matrix and $J$ the $n \times n$ matrix with all entries equal to $1,$ and let $j$ be a column $J.$ Furthermore, if $A=A(\Ga)$ then $\overline{A}=J-I-A$ is the adjacency matrix of the complement of $\Ga.$ Since $\Ga$ is regular we have $A\cdot j=k\,j$ where $k$ is the valency of $\Ga.$   If $f\in E^{*}$ then $A\cdot f=0$ implies that $f$ and $j$ are perpendicular to each other. In particular, $J\cdot f=0$ so that $\overline{A}\cdot f=(J-I-A)\cdot f=-f.$ Conversely, if $f$ is an eigenvector of $\overline{\Ga}$ with eigenvalue $-1$ then $J\cdot f=0$ and $-f=\overline{A}\cdot f=(J-I-A)\cdot f$ implies that  $A\cdot f=0.$ \dne

\medskip
\begin{thm}\label{rsu}\cite{brouwer2011spectra}
Let $\Gamma=(V,E)$ be a $k$-regular graph with vertex set $V$ and edge set $E.$ Let $\{\lambda_1,...,\lambda_n\}$ be the eigenvalues of $\Gamma.$ Then the line graph $L(\Gamma)$ is $(2k-2)$-regular graph with eigenvalues $\lambda_i+k-2$ for $1\leq i\leq n$ and $-2$ with multiplicity of $|E|-|V|.$
\end{thm}
\pf Let $\mathbb{C}$ be the field of complex numbers. Let $\mathbb{C}V$ be the vector space with basis $V$ and let  $\mathbb{C}E$ be the vector space with basis $E$ of $\Gamma.$ Thus the incident relationship  between $V$ and $E$ gives rise to the two $\mathbb{C}$-incident maps, $\varepsilon :\mathbb{C}V\longrightarrow \mathbb{C}E$ and $\sigma :\mathbb{C} E\longrightarrow\mathbb{C}V$ which are defined on the respective basis by 
$$\varepsilon (v)=\sum e  \quad \text{where}\quad e \quad \text{is an edge of}\quad v$$
and 
$$\sigma (e)=\sum v \quad \text{where}\quad v \quad \text{is endpoint of}\quad e.$$

It follows that the composition map $\nu^+=\sigma\varepsilon:\mathbb{C}V\longrightarrow\mathbb{C}V$ is given $$\nu^+(v)=\sigma\varepsilon(v)=\sigma(\sum_{e\sim v}e)=\sum_{e \sim v}\sigma(e)=k+\alpha(v)$$  for all $v\in V$ and the composition map $\nu^-=\varepsilon\sigma:\mathbb{C}E\longrightarrow\mathbb{C}E$ is given by   $$\nu^-(e)=\varepsilon\sigma (e)=\varepsilon(\sum_{e \sim v_i} v_i)= \sum_{e \sim v_i}\varepsilon(v_i)=2 e+\alpha^*(e)$$ for all $e\in E$ where $\alpha^*$ is the adjacency map of  $L(\Gamma)$. Note  by  \autoref{ma2} we have that $\nu^+$ and $\nu^-$ have the same non zero eigenvalues. Therefore from this we conclude that $\alpha^*(e)=k+\alpha(v)-2$ for some $e\in E$ and $v\in V.$\dne
\begin{cor} \label{Ali23}
Let $\Gamma$ be a $k$-regular  graph. Then $\Gamma$ is singular if and only if $k-2$ is an eigenvalue of $L(\Gamma).$
\end{cor}
\pf Suppose that $\Gamma$ is singular. Then $0$ is an eigenvalue of $\Gamma.$ Hence by \autoref{rsu} we have that $k-2$ is an eigenvalue of $L(\Gamma).$ Conversely, suppose that $k-2$ is an eigenvalue of $L(\Gamma).$ So according to \autoref{rsu} we have that $k-2=\lambda_i-2+k$ where $\lambda_i$ is an eigenvalue of $\Gamma.$   From this we deduce that $\lambda_i=0$ so that $\Gamma$ is singular. Other possibility we have that $k-2=-2$ hence $k=0$ and this gives us a contradiction.   
\dne

Next we list a few general properties of singular graphs. If $\Ga=(V,E)$ is a graph and $X\subseteq V$ then the {\it induced subgraph} $\Ga'=\Ga[X]$ is the graph $(X, E')$ where $E'$ consists of all $\{v,v'\}\in E$ with both $v$ and $v'$ in $X.$ Our first criterium is a balance condition that plays a role in applications for physical and chemical systems, see \cite{graovac1972graph} and \cite{sciriha2008nut} as a reference.

\medskip
\begin{thm}[Balance Condition]  \label{sing2} Let $\Gamma$ be a graph with vertex set $V.$ Then $\Gamma$ is singular if and only if there are disjoint non-empty subsets $X,\,Y\subseteq V$ and a function $f:X\cup Y\rightarrow \mathbb{N}$ with $f(u)\neq 0$ for all $u\in X\cup Y$ so that $$\sum_{v\sim u\in X}\,f(u)=\sum_{v\sim u\in Y}\, f(u)$$ for all $v\in V.$ In particular, if $X,Y\subset Z,$ then $\Gamma[Z]$ is singular. Furthermore, if $\Gamma$ is singular with $X,Y$ as above then  ${\rm null}(\Gamma)\geq min\{|X|,|Y|\}.$
\end{thm}

\pf Suppose that $\Gamma$ is singular of $n$ vertices. Let $h\in E_{*}$ with  $h\neq 0.$ As $A=A(\Gamma)$ is singular  over $\mathbb{Q}$ we can assume that $h_{v}$ is rational for all $v\in V,$ and after multiplying by the least common multiple of all denominators,  that $h_{v}$ is an integer for all $v.$  Let $X$ be the set of all $v$ such that $h_{v}\geq 1$ and $Y$ the set of all $v$ such that $h_{v}\leq  -1.$ Define $f^{X}$ and $f^{Y}$ in $\mathbb{C}V$ by $f^{X}_{v}=h_{v}$ for $v\in X$  and   $f^{X}_{v}=0$ otherwise, while  $f^{Y}_{v}=-h_{v}$ for $v\in Y$  and   $f^{Y}_{v}=0$ otherwise. Thus
\begin{equation}\label{sn}
A\cdot f^{X}=A\cdot f^{Y}.
\end{equation}
    For any $v\in V$ we have that  $\langle v,A\cdot f^{X}\rangle =\langle A\cdot v, f^{X}\rangle =\sum_{v\sim u\in X}\,f^{X}(u). $  Here we use that $A$ is self-adjoint, that is $$\langle h,Ak\rangle =\langle Ah,k\rangle$$ for all $h,k \in \mathbb{C}V.$ Similarly, $\langle v,A\cdot f^{Y}\rangle =\sum_{v\sim u\in X}\,f^{Y}(u).$ Hence by \autoref{sn} we have that  $\sum_{v\sim u\in X}\,f(u)=\sum_{v\sim u\in Y}\, f(u)$ for all $v\in V.$

Suppose that the above condition holds. This means   that $$\sum_{v\sim u\in X}\,f(u)=\sum_{v\sim u\in Y}\, f(u)$$ for all $v\in V.$ Now we prove that $\Gamma$ is singular. As before  $A$ is the adjacency matrix of $\Gamma.$ Note $A$ is non-singular if and only if its rows are linearly independent. Suppose that $A_v$ be the row of $A$ labelled by  a vertex $v \in V(\Gamma).$ Note we have that $$\sum_{x \in X} f(x) A_x-\sum_{y \in Y} f(y) A_y=0$$ where $f(v)\neq 0$ for all $v \in X\cup Y.$ From this we conclude that the rows of $A$ are linearly dependent and so $A$ is singular.

Now we want to prove that ${\rm null}(\Gamma)\geq min\{|X|,|Y|\}.$ Let $X, Y$ be as in the first part of the theorem and let $|X|=m$ and $|Y|=\l$ where $\l\leq m.$ As before $A$ is  the adjacency matrix of $\Gamma$ and  $A_v$ is the row of $A$ labelled by  $v$ where $v\in V.$   So by the Balance Condition we have that   $$\sum_{v\sim u\in X}\,f(u)=\sum_{v\sim u\in Y}\, f(u)$$ for all $v\in V.$  In this case $A$ has the following shape 

\[ A= \left(
\begin{array}[pos]{c}
A_{x_1}\\
A_{x_2}\\
  \vdots \\
A_{x_m}\\
A_{y_1}\\
A_{y_2}\\
  \vdots \\
  A_{y_{\l}}\\
  A_{v_{\l+m+1}}\\ 
   \vdots \\
  A_{v_{n}}\\
\end{array}
\right).
\]
 Note we do some elementary row operations for $A$ as it is shown in the following  
\[ A= \left(
\begin{array}[pos]{c}
A_{x_1}\\
A_{x_2}\\
  \vdots \\
\sum_{x_i \in X}f(x_i)A_{x_i}\\
A_{y_1}\\
A_{y_2}\\
  \vdots \\
  \sum_{y_i \in Y}f(y_i) A_{y_i}\\
    A_{v_{\l+m+1}}\\
  \vdots\\
  A_{v_{n}}\\
\end{array}
\right).
\]

By this condition and multiply $m^{th}$ row of $A$ by $-1$ and add it to $\l^{ th}$ row. In this case  we get $\l^{th}$ row equals to zero as it is shown in the following 
\[ A= \left(
\begin{array}[pos]{c}
A_{x_1}\\
A_{x_2}\\
  \vdots \\
\sum_{x_i \in X}f(x_i)A_{x_i}\\
A_{y_1}\\
A_{y_2}\\
  \vdots \\
  A_{\l}=0\\
    A_{v_{\l+m+1}}\\
  \vdots\\
  A_{v_{n}}\\
\end{array}
\right).
\]

After that we do some elementary row operations for $A$ as it is shown in the following  
\[ A= \left(
\begin{array}[pos]{c}
A_{x_1}\\
A_{x_2}\\
  \vdots \\
\sum_{x_i \in X}f(x_i)A_{x_i}\\
A_{y_1}\\
A_{y_2}\\
  \vdots \\
  \sum_{i=1}^{l-1}f(y_i)A_{y_i}\\
  A_{\l}=0\\
    A_{v_{\l+m+1}}\\
  \vdots\\ 
  A_{v_{n}}\\
\end{array}
\right).
\]

Then we subtract some of $f(x_i)A_{x_i}$ from $m^{th}$ row and add these   to the $(\l-1)^{th}$ row so that the Balance condition hold. After that we multiply  $m^{th}$ row by $-1$ and add it to $(\l-1)^{th}$ row we get another row equals to zero. Keep doing this list of row elementary row operations we get $\l^{th}$ rows of $A$ equal to zero. Hence from the above we conclude that ${\rm null}(A)\geq \l=|Y|.$  \dne

\textit{ Example 1}: Let $\Gamma=(V,E)$ be a graph. Suppose that  $w,u \in V$ such that $w \nsim u,$ and $w$ and $u$ have the same neighbour set. In this case put $X=\{w\}, Y=\{u\}$ and $f(u)=f(w)=1$ while $f(v)=0$ for $u\neq v\neq w.$ Then $f$ has the property of the theorem. More directly of course, $\alpha(w)=\alpha(u)$ implies that $0\neq w-u \in E_*.$

\textit{ Example 2}: Let $\Gamma=C^{n}$ be an $n$-cycle on $V=\{1,2,...,n\}.$ Then it is easy to compute the eigenvalues of $\Gamma.$ These are the numbers $\lambda_r=2cos(\frac{2\pi r}{n})$ where $r=0,1,2,...,n-1,$ see \cite{borovicanin2009nullity} as a reference. In particular, $C^{n}$ is singular if and only if $n$ is divisible by $4.$ If $4$ does divide $n$ we may take $X=\{a\in V\,\,: \,\,a\equiv 0 {\rm \,\,or\,\,} 1\pmod 4\},$ $Y=\{b\in V\,\,: \,\,b\equiv 2 {\rm \,\,or\,\,} 3\pmod 4\}$ and $f(v)=1$ for all $v\in V.$

\textit{ Example 3}: Let $\Gamma=P^{n}$ be a path on $n$ vertices, $V=\{1,2,...,n\}.$ Then it is easy to compute the eigenvalues of $\Gamma.$ These are the numbers $\lambda_r=2cos(\frac{\pi r}{n+1})$ where $r=1,2,...,n.$ In particular, $P^{n}$ is singular if and only if $n$ is odd see \cite{borovicanin2009nullity} as a reference.

\section{\sc Cayley Graphs Spectra}
Let $G$ be a finite group and let $H$ be a subset of $G$ with the following properties:
\begin{enumerate}[\qquad(i)]
\item $H$ generates $G,$ 
\item $H=H^{-1}$ and 
\item $1\not\in H.$ 
\end{enumerate}
Then $H$ is called a {\it connecting set} in $G.$
Define the Cayley graph $\Ga={\rm Cay}(G,H)$ with vertex set $V(\Ga)=G=V$ so that $v\in V$ is adjacent to $v'$ if and only if there is some $h$ in $H$ for which $hv'=v.$ Note, if $hv'=v$ then $v'=h^{-1}v$ with $h^{-1}\in H$ so that $v'$ is adjacent to $v.$ In particular, $H$ is the set of vertices which are adjacent to $1\in G$ and the conditions (i)-(iii) imply that $\Ga$ is a connected, undirected and loop-less graph with vertex set $G.$

The multiplication of vertices  on the right by an element $g$ in $G,$ that is $v\mapsto vg$ for $v$ in $V,$ induces an automorphism on $\Ga.$ To see this let $v,v'\in V.$ If $v\sim v'$ then $v'=h^{-1}v$ for some $h\in H$ and so $v'g=h^{-1}(vg)$ giving that $vg\sim v'g.$ And conversely, if $v\not \sim v'$ then $vg\not \sim v'g.$
This automorphism gives rise  to the {\it right regular representation} $\rho_{r}\!:\,G\to {\rm GL}(\cp V)$  of $G$ given by $\rho_{r}(g)(v)=vg.$

In $\Ga={\rm Cay}(G,H)$ the vertices $v'$ adjacent to a given $v$ are of the shape $v'=h^{-1}v$ with $h\in H.$ Therefore $$(*)\quad \ga(v)=\sum_{h\in H}\,h^{-1}v$$ and in particular, $\rho_{r}(g)$ commutes with $\alpha$ for all $g\in G.$ 

However, we can understand $\ga$ also in terms of the {\it left regular representation} $\rho_{\ell}\!:\,G\to {\rm GL}(\cp V)$  of $G$ given by $\rho_{\ell}(g)(v)=g^{-1}v.$  It follows that   \begin{equation}\label{eigen1}  
\alpha=\sum_{h\in H}\,\rho_l(h)\,\,\, \text{as a map }\mathbb{C}V\to \mathbb{C}V.
 \end{equation} 

It can be seen easily that the multiplication of vertices on the left, that is $v\mapsto g^{-1}v$ for $v$ in $V,$ is an automorphism of $\Ga$ if and only if $g^{-1}H=Hg^{-1}.$ Therefore $\rho_{\ell}(g)$ commutes with $\alpha$ if and only if $g^{-1}H=Hg^{-1},$ or equivalently, $gH=Hg.$
We say that $H$ is {\it normal,} and that $\Ga={\rm Cay}(G,H)$ is a {\it normal} Cayley graph, if $gH=Hg$ for all $g\in G.$

The right regular  action is transitive on vertices and only the identity element fixes any vertex. This is therefore the regular action of $G$ on itself. This property characterises Cayley graphs.

\begin{thm}\label{zu5}(Sabidussi's Theorem)\cite[p.~48]{godsil2013algebraic}
 Let $\Gamma=(V,E)$ be a graph. Then $\Gamma$ is a Cayley graph if and only if $Aut(\Gamma)$ contains a subgroup $G$ which is regular on $V.$
\end{thm}
Now we list results that deal with the spectrum of Cayley graphs:

\begin{thm} \label{noor41}
Let $G$ be a finite group and let  $\rho_1,\cdots, \rho_s$ be the  set of all inequivalent irreducible representations of $G.$  Then $\lambda$ is an eigenvalue of $Cay(G,H)$ if and only if there is some $\rho_i$ such that  $\sum_{h\in H} \rho_i(h)-\lambda$ is singular. 
\end{thm}

\pf Let $U_1,...,U_s$ be the  irreducible $G$-modules. Let $E_1,...,E_t$ be the eigenspaces of $\alpha.$ By \autoref{noor47}   we have that $E_1,...,E_t$ are $G$-invariant (under multiplication on the right) and so each $E_j$ can be decomposed into $$E_j=m_{j1}U_1\oplus...\oplus m_{js}U_s,$$ as before. Now $E_j$ is the eigenspace of $\alpha$ for the eigenvalue $\lambda$ if and only if $\alpha-\lambda$ is singular on $E_j.$ This in turn implies that $\sum_{h \in H}\rho_l(h)-\lambda$ is singular on $E_j.$ Let $U_i$ be an irreducible $G$-module that appears in $E_j.$  Then $\sum_{h \in H}\rho_i(h)-\lambda$ is singular. 

Conversely, if $\sum_{h \in H}\rho_i(h)-\lambda$ is singular on $U_i$ then $U_i$ appears in the decomposition of $$\mathbb{C}G=E_1\oplus...\oplus E_t$$ as this the regular $G$-module, and so $\alpha-\lambda$ is singular.\dne

\begin{thm}\cite{krebs2011expander, diaconis1981generating, brouwer2011spectra}  \label{noor14}
Let $U$ be an irreducible sub-module of $\mathbb{C}V=\mathbb{C}G,$ (by right multiplication).  Suppose $H$ is a normal connecting set of $G.$ Let $\Gamma=Cay(G,H).$ Then $\alpha(U)= U$ and furthermore $U$ is contained in the eigenspace of $\alpha$ for $$\lambda=\frac{1}{\chi(1_G)} \sum_{h \in H}\chi(h)$$ where $\chi$ is the irreducible character corresponding to $U.$
\end{thm}
\pf Let $\mathbb{C}V$ be the vertex $G$-module of $\Gamma.$   Let $\rho=\rho_i: G\rightarrow GL(U).$  Then  by using \autoref{eigen1} we have  that  $\alpha=\sum_{h\in H} \rho_l(h)$ where   $\rho_l(h)$ is the left regular representation of $G.$  Thus we compute
\begin{eqnarray}
\alpha \circ \rho(g)&=&\sum_{h\in H} \rho_l(h)\,\rho(g)\nonumber\\
&=& \sum_{h\in H}\sum_{j=1}^s m_j \rho_j(h)\rho(g)\nonumber\\
&=&\sum_{h\in H}\rho_i(hg)\nonumber\\
&=&\sum_{h\in H} \rho_i(g h g^{-1} g)\nonumber\\
&=&\sum_{h\in H} \rho_i(gh)\nonumber\\
&=&\rho_i(g)\sum_{h\in H} \rho_i(h) \nonumber\\
&=&\rho_i(g)\sum_{h\in H}\sum_{j=1}^s m_j\rho_j(h) \nonumber\\
&=&\rho(g) \circ \alpha \nonumber
\end{eqnarray}
for all $g\in G$. 

Therefore by using  \autoref{noor4} we have that $\alpha(u)=\lambda u $ for some $\lambda\in \mathbb{C}$ and all $u \in U.$ So we have that $U \subseteq E_i$ for some $i$ where $E_i$ is the eigenspace corresponding to the eigenvalue $\lambda.$ Note we have that  $$m \lambda={\rm tra}(\lambda u)={\rm tra}(\alpha(u))={\rm tra}(\sum_{h\in H} \rho(h)(u))=\sum_{h\in H} \chi(h)$$ where $m=\chi(1_G).$  So we have that $$\lambda=\frac{1}{\chi(1_G)} \sum_{h \in H} \chi(h).$$ For the multiplicity  each irreducible character $\chi$ there are $\chi(1_G)$  copies of $U$ in $\mathbb{C}V$ and on each copy $\alpha$ acts as $\lambda id_{U}.$ Therefore $\lambda$ has multiplicity  $(\chi(1))^2.$ \dne 

\subsection{ Singular Cayley graphs}

 As before  $G$ is  a finite group and $H$ a connecting set of $G.$ Let $\Gamma$ be the  graph $\Gamma=Cay(G,H).$ We denote  by $K$ a subgroup of $G.$ Note, we say that $H$ is \textit{vanishing}  on the irreducible character $\chi$ if $\sum_{h\in H}\chi(h)=0.$  In this section we investigate conditions for $\Gamma$ to be singular.  
 
  \begin{thm}\cite{klotz2010integral}\label{sincay3}
 Let $H$ be a connecting set of a group $G$ and let $H$ be a union of left cosets of a  non-trivial subgroup $K$ in $G.$ Suppose there is some element $k \in K$ and $1$-dimensional character $\chi$ of $G$ such that $\chi(k)\neq 1.$ Then we have that $\sum_{h \in H}\chi(h)=0.$ In particular, $\Gamma=Cay(G,H)$ is singular.    
 \end{thm}
 This theorem is mentioned in \cite{klotz2010integral} for a subgroup of an additive abelain group. However, the following proof is my version for a union of cosets of a non-trivial subgroup of the  group in general.

  \pf  Suppose that $H=a_1K\cup a_2K\cup...$  for $a_1,a_2,... \in G$ and $\chi$ is an $1$-dimensional character of $G$ with $\chi(k)\neq 1$ for some $k \in K.$  Then we have that  
\begin{eqnarray}
\sum_{h\in H}\chi(h)&=&\sum_{a_i\in G}\chi(a_iK)\nonumber\\
&=&\sum_{a_i\in G}\chi(a_iKk)\nonumber\\
&=& \sum_{a_i\in G}\chi(a_iK)\chi(k)\nonumber\\
&=&\sum_{h\in H}\chi(h)\chi(k).\nonumber
\end{eqnarray}
Hence we have that $\sum_{h\in H}\chi(h)(1-\chi(k))=0.$ Since $\chi(k)\neq 1$ so that $\sum_{h\in H}\chi(h)=0.$\dne
  \begin{prop}\label{myp1}
Let $G$ be a  group with normal subgroup $K$ and a homomorphism $$\varphi:G\rightarrow G/K.$$ Suppose that $H$ is a subset of $G$ such that:

(1) $\varphi(H)$ is vanishing in $G/K$ for some character $\chi$ of $G/K.$

(2) There is a constant $c$ such that every coset of $K$ in $G$ meets $H$ in $0$  or $c$ elements. 

Then $H$ is vanishing in $G.$
\end{prop}
\pf Let $\chi$ be the irreducible character of $G/K$ on which $\chi$ is vanishing.  Then we have that \begin{equation}\label{ho_1}
0=\sum_{h \in H} \chi(\varphi(h))=\sum_{g_iK \cap H\neq \phi}\chi(g_iK)
\end{equation} 
where $g_1K \cupdot g_2K \cupdot ...\cupdot g_mK=G.$ Hence by \autoref{ho_1} we have that $$\sum_{h\in H}\tilde{\chi}(h)=c\sum_{g_iK\cap H\neq \phi}\chi(g_iK)=0$$ where $\tilde{\chi}$ is the lift character corresponding to $\chi.$\dne

\begin{thm}\label{noor15}
Let $G$ be a finite group and let $H$ be a connecting set of $G.$ Let $\rho_1,..,\rho_s$ denote the irreducible representations of  $G.$ Then $\Gamma=Cay(G,H)$ is singular if and only if there exists some $i$  such that $\sum_{h\in H}  \rho_i(h)$  is singular.   
\end{thm}

\begin{thm}\label{noor16}
Let $G$ be a finite group and let $H$ be a connecting set and normal subset of $G.$ Then $Cay(G,H)$ is singular if and only if  there is an irreducible character $\chi$ of $G$ such that $\sum_{h \in H} \chi(h)=0.$ In particular, we have that ${\rm null}(\Gamma)\geq (\chi(1_G))^2.$   
\end{thm}

\pf Suppose that $\Gamma$ is singular. Then we have that $0$ is an eigenvalue of $\Gamma.$ Note $H$ is normal subset of $G$ and so by  \autoref{noor14} each eigenvalue of $Cay(G,H)$ is given by  $$\lambda=\frac{1}{\chi(1_G)}\sum_{h\in H}\chi(h)$$ where $\chi$ is an irreducible character of $G.$  Hence  we have that $\sum_{h\in H}\chi(h)=0$ for some irreducible characters of $G.$ 

Suppose that $\sum_{h\in H}\chi(h)=0$ for some irreducible character of $G.$  Then by  \autoref{noor14} we have that $\lambda=0$ for some eigenvalues of   $Cay(G,H)$ so $Cay(G,H)$ is singular. Furthermore, we have that $E_*$ contains the module of $\chi$ with multiplicity of $(\chi(1_G))^2.$
\dne

\begin{cor}\label{cayly3}
Suppose  $G$ is non-abelian  simple group. Suppose $H$ is any subset of $G$ with $1_G \not\in G, H=H^{-1}$ and $H$ is normal. Then nullity of $\Gamma=Cay(G,H)$ is either $0$ or $\geq m^2$ where $m \neq 1$ is the least degree of an irreducible character of $G.$  
\end{cor}

 \begin{thm}\label{sincay2}
Let  $H$ be a connecting set in the  group $G$ and suppose that $H$ is a union of right cosets of the subgroup $K$ of $G$ with $|K|\neq 1.$ Then $A(\Gamma)$ is of the form $A(\Gamma^*)\otimes J$ where $\Gamma=Cay(G,H),$  $\Gamma^*$ some graph defined on the right cosets of $K$ in $G$ and $J$ is the $|K| \times |K|$ matrix with all entries equal to $1.$
\end{thm}
{\sc Comments:}(1) If $H$ is a union of left  cosets of $K$  then it is a union of right cosets since $H=H^{-1}.$ However, $aK\cup Ka^{-1}$ may  not be a union of left or right  cosets. 

(2) Note $\Gamma^*$ is a coset graph and in general may not be a Cayley graph. In some cases $\Gamma^*$ is a Cayley graph, for instance if $K$ is normal.
 
\textit{Proof}  of \autoref{sincay2}: Suppose that $H=Ka_1\cup Ka_2\cup Ka_3\cup...$ for some $a_1, a_2,... \in G.$ We want to prove that any two elements in the same right coset of $K$ are not adjacent and if an element in $Kg_i$ is adjacent to an element in $Kg_j$ then all elements in $Kg_i$ are adjacent to all  elements in $Kg_j.$
 
 Let $x , \tilde{x} \in Kg_i.$ Suppose that $x \sim \tilde{x}.$ Hence by the  Cayley graph definition we have $ \tilde{x}=h^{-1} x$ for some $h \in H.$ So we have that $\tilde{k_1}g_i=a'k_2 k_1 g_i$ for some $a' \in G$   so that $\tilde{k_1}=a'k_3$ for some $\tilde{k_1}, k_1, k_2, k_3 \in K.$ This give us a contradiction as the right cosets of $K$ are disjoint. From the above we conclude that the elements of the same right cosets are non-adjacent.
 
 Now let $x, \tilde{x} \in Kg_i$ and $y, \tilde{y} \in Kg_j.$ Suppose that $x \sim y$ in $\Gamma$ and  we want to show that $\tilde{x} \sim \tilde{y}$ in $\Gamma.$ Note by the Cayley graph definition we have that $y=h^{-1} x$ for some $h \in H.$ In this case we have that $k_2g_j=a'k_3 k_1g_i$ so $g_j=k_2^{-1}a'k_4g_i.$ So we have that 
 \begin{eqnarray}
  \tilde{y} &=&\tilde{k_2}g_j\nonumber\\
  &=&\tilde{k_2} k_2^{-1}a'k_4g_i\nonumber\\
  &=&\tilde{k_3}a'k_4g_i
 \end{eqnarray}
for $k_1, k_2, k_3, k_4, \tilde{k_2}, \tilde{k_3} \in K.$ Now assume that $\tilde{x}=k_4g_i$ hence $\tilde{x} \sim \tilde{y}.$ From the above we conclude  that all  elements in $Kg_i$ are adjacent to all elements in $Kg_j.$ From the above we deduce that $\Gamma$ is imprimitive graph and $$A(\Gamma)=A(\Gamma^*)\otimes J.$$ \dne

\begin{cor}\label{sincay1}
If $H$ is a connecting set in the  group $G$ and if $H$ is a union of right cosets of the subgroup $K \subseteq G$ with $|K|\neq 1,$ then  $\Gamma=Cay(G,H)$ is singular and ${\rm null}(\Gamma)\geq \frac{|G|}{|K|}\cdot(|K|-1).$
\end{cor}
\pf Note $J$ is singular and the eigenvalues of $J$ are $|K|$ with multiplicity of $1$ and $0$ with multiplicity of $|K|-1.$ Hence we have that $$Spec(\Gamma)=Spec(\Gamma^*)\otimes (|K|,\underbrace{0,0,...,0)}_{|K|-1}.$$ From this we conclude that ${\rm null}(\Gamma)\geq  \frac{|G|}{|K|}\cdot(|K|-1).$\dne

\subsection{Singular Cayley Graphs over a Cyclic Group}

In this section we derive  simple conditions which characterise   singular Cayley graphs over a cyclic group. Note that a Cayley graph over a cyclic group  is also called a \textit{circulant graph}. Let $C_n=\langle a\rangle$ be a cyclic group of order $n$ and let $H$ be a connecting set of $C_n.$ Denote the Cayley graph $Cay(C_n,H)$ by $\Gamma.$ It is clear that $H$ is a normal subset of $C_n.$ 

Note each irreducible representation of $C_n$ has degree $1$ see \cite{james2001representations} as a reference. Then $\rho_i(a)=\omega^{i-1}$ for $i=1,2,..,n$  are the complete list of the irreducible representations of $C_n$ and the same time these are the  irreducible characters of $C_n$ where $\omega$ is a primitive $n^{th}$ root of unity.  Therefore by \autoref{noor14} we deduce that each eigenvalue of $\Gamma$ is a certain sum of $n^{th}$ roots of unity. Note that in this work the irreducible character $\chi_i$ for $1\leq i\leq n$ of $C_n$ is generated by $\omega^{i-1}$ where $\omega$ is a fixed primitive $n^{th}$ root of unity  and the corresponding eigenvalue of $\Gamma$ will be $\lambda_i.$

  Let $\l$ be a positive  integer and let $\Omega_{\l}$ be the group of  $\l^{ th}$ \textit{roots of unity}, that is  $\Omega_{\l}=\{z\in \mathbb{C}\backslash\{0\} :z^{\l}=1\}.$ Then $\Omega_{\l}$ is a cyclic group of order $\l$ with  generator  $e^\frac{2\pi i}{\l}.$ Note this is not  the only generator of $\Omega_{\l},$  indeed any power $e^\frac{2\pi im}{\l}$ where $gcd(\l,m)=1$ is a generator too. A generator of $\Omega_{\l}$ is called \textit{a primitive $\l ^{th}$ root of unity}. 
  
     Let $n$ be a positive  integer and let $\Phi_n(x)$ denote the $n^{th}$ \textit{cyclotomic polynomial}. Then $\Phi_n(x)$ is the unique irreducible integer polynomial with leading  coefficient $1$ so that $\Phi_n(x)$ divides  $x^n-1$ but  does not divide  of $x^k-1$ for any $k<n$. Its roots are all primitive $n^{th}$  roots of unity. So  $$\Phi_n(x)=\sum_{1\leq m<n}(x-e^\frac{2\pi im}{n})$$ where $gcd(m,n)=1.$ 
     
     \begin{lem}\cite[Lemma 3.1.1]{szabo2004topics}\label{gh37}
If $n$ is a prime power, $n=p^m,$ if  $\omega$ is a primitive $n^{th}$ root of unity and if $a(1),...,a(k)$ are integers with
\begin{equation}\label{cay6}
\omega^{a(1)}+...+\omega^{a(k)}=0
\end{equation}
then $k$ is  a multiple of $p.$
 
\end{lem}

\begin{lem}\cite[Lemma 3.1.3]{szabo2004topics}\label{cycl}
Let $\omega$ be a primitive $n^{th}$ root of unity and let $a(1),...,a(k)$ be integers. If $n$ is a product of two prime powers, say $n=p^e q^f$, and if
$$\omega^{a(1)}+...+\omega^{a(k)}=0,$$ then $$\omega^{a(1)}+...+\omega^{a(k)}=\l (1+\delta+...+\delta^{p-1})+r(1+\varepsilon+...+\varepsilon^{q-1}),$$ where $\l, r$ are sums of powers of $\omega,$ and $\delta, \varepsilon$ are primitive $p^{th}$ and $q^{th}$ roots of unity respectively. 
\end{lem} 

As a consequence to \autoref{gh37} and \autoref{cycl} we have the following result:
   \begin{thm}\label{gh35}
 Let $C_n=\langle a\rangle $ and $H$ a connecting  set of $C_n.$  Let $\Gamma$ be the graph $\Gamma=Cay(C_n,H)$ and let $\Omega_n$ be the  group of  $n^{th}$ roots of unity. For $i=1,2,3,...,n$ consider the  homomorphism  $$\varphi_i:C_n\rightarrow \Omega_n$$ given by $\varphi_i(a^m)=\omega^{(i-1)m}$  where $\omega$ is a primitive $n^{th} $ root of unity and $0\leq m\leq n-1.$ Then $\Gamma$ is singular graph if  the multi-set $$\varphi_i(H)=\{\varphi_i(h_1),...,\varphi_i(h_k)\},$$ with $|H|=k,$  is a union  of  cosets of some non-trivial subgroup $\Upsilon\subseteq \Omega_n$ for some $i.$
 \end{thm}

\begin{cor}\label{cycl1}
Let $G$ be a  group with normal subgroup $K$ such that $G/K$ is abelian. Let $H$ be a connecting set of $G$ and let $\Gamma=Cay(G,H).$ Suppose that every coset of $K$ in $G$ meets $H$ in  exactly $c$ elements for some $c.$ Then $\Gamma$ is singular with nullity $\geq |G/K|-1.$
\end{cor}
\pf We need to show that $A\cong \sum_{h \in H}\rho_{\l}(h)$ is singular and of nullity $\geq |G/K|-1,$ where $A$ is the adjacency matrix of $\Gamma$ and $\rho_{\l}$ is the left regular representation of $G.$ Then  we can decompose  $$\sum_{h \in H} \rho_{\l} (h)=m_1\sum_{h \in H} \rho_1(h)\oplus...\oplus m_s\sum_{h \in H}\rho_s(h)$$ where $\rho_1,...,\rho_s$ are the irreducible representations of $G$ and $m_1,...,m_s$ are their degrees respectively. For this it is sufficient to show that $\sum_{h \in H}\rho_i(h)$ is singular for some $i.$ Now let $\chi_1,...,\chi_{|G/K|}$ be the irreducible characters of $G/K.$ Then we have that $\sum\chi_j(gK)=0$ for all non-trivial irreducible characters of $G/K$ where the sum over all the cosets of $K$ in $G.$ Hence we have that $\sum_{h \in H}\tilde{\chi}_j(h)=c\sum\chi_j(gK)=0$ where $\tilde{\chi}_j$ is the lift character corresponding to $\chi_j.$ From this we conclude that $\sum_{h\in H}\rho_j(h)$ is singular where $\rho_j$ is the irreducible representation of $G$ which is corresponding to $\tilde{\chi}_j.$ So  by \autoref{noor41} we have that $A$ is singular with nullity $\geq |G/K|-1$ as there are $|G/K|-1$ non-trivial character for $G/K.$\dne

As before  $C_n=\langle a \rangle$ is a cyclic group of order $n$ and  $H$  a connecting set of $C_n.$ Let $\Gamma$ be the graph $\Gamma=Cay(C_n,H)$ and  let $H^*$ be the set of all $0<m\leq n-1$ such that $H=\{a^m:m\in H^*\}.$ Now consider the polynomial $$\Psi_\Gamma(x)=\sum_{m\in H^*}x^m$$ associated to $\Gamma.$ Note that $\Psi_\Gamma$ depends on the choice of the generator $a.$ If $a'$ is some other generator, then $a=(a')^r$ for some $r$ with $gcd(r,n)=1.$ Therefore $(H')^*\subseteq \{1,2,...,n-1\}$ given by $(H')^*\equiv r H^* \mod n$ and so $$\Psi'_\Gamma(x)=x^r\sum_{m\in H^*}x^m\equiv x^r \Psi_\Gamma(x) \mod x^n.$$ 
   \begin{thm}\cite{lal2011non, kra2012circulant}\label{gh43}
Let  $C_n=\langle a\rangle $ be a cyclic group of order $n$ and let $\Gamma=Cay(C_n,H)$ be the Cayley graph for the connecting set $H\subset C_n.$  Let $\Psi_\Gamma(x)$ be  the polynomial associated to  $\Gamma$ for some generator of $C_n.$ Then $\Gamma$  is singular if and only if $\Phi_d(x)$ divides $\Psi_\Gamma(x)$ for some divisor $d$ of $n$ with $1<d\leq n$ where $\Phi_d(x)$ is the $d^{th}$ cyclotomic polynomial. 
Furthermore, let $d_1,d_2,...,d_l$ be the divisors of $n.$ Then  we have that ${\rm null}(\Gamma)=\sum \varphi(d_j)$ where the sum is over all $d_j$ such that $\Phi_{d_j}(x)$ divides $\Psi_\Gamma(x).$
   \end{thm}
   Note, this proof in our version based on our techniques. 
   
       \pf Let $\lambda=\lambda_i$ be any eigenvalue of $\Gamma.$ Thus by \autoref{sincay3} we have that
       \begin{eqnarray}\label{my1}
       \lambda_i &=& \sum_{h\in H}\chi_i(h)\nonumber\\
       &=& \sum_{m\in H^*}\chi_i(a^m)\nonumber\\
       &=& \sum_{m\in H^*}\chi_i(a)^m\nonumber\\
       &=& \sum_{m\in H^*}(\omega^{i-1})^m\nonumber\\
       &=& \Psi_\Gamma(\omega^{i-1})\nonumber
       \end{eqnarray}
       
        where  $\omega$ is a primitive  $n^{th}$ root of unity. Now if  $\lambda_i=0,$ then  $\Psi_\Gamma(\omega^{i-1})=0$ and so; if $\omega^{i-1}$ is  a primitive  $n^{th}$ root of unity then we have that  $\Phi_n(x)|\Psi_\Gamma(x)$ as $\Psi_\Gamma(x)$ and $\Phi_n(x)$ have a common  root and if  $\omega^{i-1}$ is not a primitive $n^{th}$ root of unity.   In this case  $\omega^{i-1}$ is a primitive $r^{th}$  root of unity for some divisor $r$ of $n$ where $1<r<n.$ Hence   we have that $\Phi_r(x)|\Psi_\Gamma(x)$ as   $\Psi_\Gamma(x)$ and $\Phi_r(x)$ have a common  root. Therefore  for both cases we have that $\Phi_d(x)$ divides $\Psi_\Gamma(x)$ for some divisor $d$ of $n$ with $1<d\leq n.$

     Conversely, suppose $\Phi_d(x)$ divides $\Psi_\Gamma(x)$ for some divisor $d$ of $n.$  Then we have that $\Phi_d(\omega^*)=0 $ where $\omega^*$ is a primitive $d^{th}$ root of unity. So we have that $\Psi_\Gamma(\omega^*)=0$ then  $\lambda_i=0$ for some $i.$ By this we deduce that $\Gamma$ is singular.

 By the second part of the proof we have that $\lambda_i=\Psi_\Gamma(\omega^*)=0$ if and only if $\omega^*$ is a primitive $d^{th}$ root of unity for some divisor $d$ of $n.$ In this case we have that $\varphi(d)$ of primitives $d^{th}$ root of unity. Hence we deduce that  ${\rm null}(\Gamma)=\sum\varphi(d_j)$ where the sum is over all divisors of $n$ such that $\Phi_{d_j}(x)$ divides $\Psi_\Gamma(x).$ \dne
 
 \begin{thm}\label{gh7}
Let  $\Gamma$ be a vertex transitive graph on $p$ vertices with at least one edge where $p$ is a prime number. Then $\Gamma$ is non-singular.
\end{thm}

\pf Let $V$ be the vertex set of $\Gamma,$ with  $|V| =p$ and $p$ is a prime number. Let $G$ be a vertex transitive group on $\Gamma.$ By \autoref{zu6} we have that $$|G|=|V| \cdot |G_v|$$ for some $v \in V.$ So by Sylow's Theorem there exist a subgroup  $K$  of $G$ with $|K|=p.$  Note $K$ is cyclic. Now apply \autoref{zu6} again we have that $|K|=|v^k|.|K_v|$ for $k \in K.$ Note we have that $|K_v|= 1_G,$ so $K$  acts regularly on $V.$ Hence by Sabidussi's \autoref{zu5} we have that $\Gamma$ is a Cayley graph $Cay(K,H)$ for some connecting set $H$ of $K.$    Suppose   for  contradiction that $\Gamma$ is singular. As $p$ is a prime number, according to \autoref{gh43} we have that $\Phi_p(x)$ divides $\Psi_\Gamma(x).$ So   there is $Q(x)\in \mathbb{Q}[x]$ such that $$\Psi_\Gamma(x)=\Phi_p(x)\times Q(x).$$ This  gives us a contradiction  as $\Phi_p(x)$ has degree  $p-1$  and $Q(x)\neq 0$ but the maximum  degree of $\Psi_\Gamma(x)$ is less than $p.$
\dne
  
\section{ Vertex Transitive Graphs}

Next we consider the singularity of a vertex transitive graph.  A graph $\Gamma$ is said to be \textit{vertex transitive } if its automorphism group acts transitively on its vertex set. In other words, for any two vertices $u,v$ of $\Gamma$ there is $g\in Aut(\Gamma)$ such that $v^g=u.$  It is clear that vertex transitive graphs are regular. In this section we compute  the spectrum of a vertex transitive graph in  terms of  the irreducible characters of a transitive group of automorphisms.

As before $\Gamma$  is a simple connected graph with vertex set $V.$ Let $G$ be a vertex transitive group of automorphisms of $\Gamma.$ Let $U_{1},...,\,U_{s}$ be  the irreducible modules of $G$ with corresponding characters $\chi_{1},...,\,\chi_{s}.$ Let $E_{1},...,\,E_{t}$ be the eigenspaces of $\alpha$ with  corresponding eigenvalues $\lambda_{1},...,\,\lambda_{t}.$ Let $m_{j,i} $ be the multiplicity of $U_i$ in $E_j.$

  \begin{thm}\label{zu12} 
 Let $G$ be a group of automorphisms of the graph $\Gamma$ which acts transitively on $V=V(\Gamma).$   Consider $\tau(g):={\rm tra}(g\alpha)$ for $g\in G$ where $\alpha$ is the adjacency map of $\Gamma.$   Then $\tau(g)$ is a class function and $\langle\tau,\chi_{i}\rangle=\sum_{j=1}^{t} \,m_{j,i}\lambda _{j}.$  If the permutation action of $G$ on $V$ is multiplicity-free then the following hold 
  
(i) Every  eigenvalue of $\Gamma$  is of the shape  $\langle\tau,\chi_{i}\rangle$ for some $i$ with multiplicity of $\chi_i(1_G).$ 

(ii)  $\Gamma$ is singular if and only if $\sum\,\chi_{i}(1)<|V|$ where the sum runs over all characters $\chi_i$  with $\langle \tau,\chi_{i}\rangle\neq 0.$
\end{thm} 

{\sc Comments:} 1. It is clear that $\tau(g)={\rm tra}(g\alpha)$ is the number of times a vertex $v \in V$ is adjacent to its image $v^g$ under $g.$ In particular, $\tau(1_G)=0.$

(2) A permutation character $\psi$ of $G$ is \textit{multiplicity-free} if and only if each irreducible character of $G$ appears with multiplicity $\leq 1$ in $\psi.$  In particular,  if $\psi$ is multiplicity-free  then $G$ is transitive.

(3) As before $G,$ is transitive on the vertex set $V$  of $\Gamma.$  Therefore the permutation   representation of $G$ on $V$ is a sub-representation of the regular representation. Therefore $\sum_{j=1...t} m_{j,i}\leq \dim(U_{i})$ for all $i.$  For instance, if $G$ is abelian then all non-zero eigenvalues are of the form $\lambda _{j}=\langle \tau,\chi_{i}\rangle $ for some $i.$

(4) If $\Gamma=Cay(G,H)$ where $H$ is a normal connecting set of $G$ then $\langle \tau, \chi_i\rangle= m_i\lambda_i$ where $\lambda_i$ is an eigenvalue of $\Gamma$ and  $m_i$ is the dimension of $U_i.$  Here the  multiplicity of $\lambda_i$ is $m_i^2.$ Since by \autoref{noor14} we have that  $U_i$ appears in $\mathbb{C}V$ with multiplicity $m_i$ and so on $m_iU_i$ appears in one eigenspace.

(5) We consider the case where   the permutation character $\psi$ of $G$ is not multiplicity free.  Let $r$ be the multiplicity of $\chi_i$ in $\psi.$ Define $\tau^l(g)={\rm tra}(\alpha^lg)$ for some $l \in \mathbb{N}.$ Then, using the same ideas as in the proof of the theorem, we have 
\begin{eqnarray}\label{ver1}
\sum_{j=1}^{t} \,m_{j,i}\lambda _{j}&=&\langle\tau,\chi_{i}\rangle \nonumber\\
\sum_{j=1}^{t} \,m_{j,i}\lambda _{j}^2&=& \langle\tau^2,\chi_{i}\rangle\nonumber\\
\vdots \nonumber\\
\sum_{j=1}^{t} \,m_{j,i}\lambda _{j}^r&=&\langle\tau^r,\chi_{i}\rangle\nonumber\\
\end{eqnarray}
where $r\leq m_i.$ These are additional equations to determine the spectrum of $\Gamma.$ Please see the example of the Petersen graph with the General Affine Group.

\textit{Proof} of \autoref{zu12} : As before    let $\mathbb{C}V$ be the vertex module of $\Gamma$ and $\alpha$ be the adjacency map of $\Gamma.$ 

First we show that $\tau(g)$ is a class function. Note by \autoref{zu} we have that $\alpha h^{-1}gh= h^{-1}\alpha gh$ for $h\in G.$  Since 
\\[-60 pt]

\begin{eqnarray}
{\rm tra}(\alpha h^{-1}gh)&=&{\rm tra}(h^{-1}\alpha gh )\nonumber\\
&=&{\rm tra}(\alpha g hh^{-1})\nonumber\\
&=& {\rm tra}(\alpha g)\nonumber
\end{eqnarray}

 we have that  $\tau(g)={\rm tra}(\alpha g)$ is a class function. So we can write $\tau$  in the following shape 
\\[-40 pt]

\begin{equation}\label{rs}
\tau(g)=\langle \tau(g),\chi_1(g)\rangle \chi_1(g)+...+\langle \tau(g),\chi_s(g)\rangle \chi_s(g)
\end{equation} 

 as $\chi_1,...,\chi_s$ is an orthonormal  basis of the vector space  of all class functions.

 Let $\pi_{1},...,\pi_{t}$ be the projections  $\pi_{j}\!:\,\mathbb{C} V\to \mathbb{C}V$ with $\pi_j(\mathbb{C}V)\subseteq E_j.$ Since $G$ preserves eigenspaces and commutes with the $\pi_j$ (in both cases as $G$ commutes with $\alpha)$ we have $g\alpha\pi_{j}=g\lambda_{j}\pi_{j}=\lambda_{j}g\pi_{j}.$ Since $\pi_{1}+...+\pi_{t}= id$ we have \begin{eqnarray}
\langle\tau, \chi_{i}\rangle &=&\langle {\rm tra}(\alpha g),\chi_i \rangle\nonumber\\
&=& \langle\sum_{j}{\rm tra}(g\alpha\pi_{j}),\chi_{i}(g)\rangle\nonumber\\
&=&\sum_{j}\lambda_{j}\langle{\rm tra} \,g\pi_{j},\chi_{i}(g)\rangle\nonumber\\
&=&\sum_{j}\lambda_{j}\,m_{ji}.\nonumber
\end{eqnarray} Note, if the permutation character of $G$ on vertices is multiplicity-free then $0\leq m_{ji}\leq 1$ for all $j,i$ and for every $i$ there is at most one $j$ with $m_{ji}=1.$ Hence $\lambda_{j}=\langle\tau,\chi_{i}\rangle$ for such a pair. By the same argument, $\lambda_{j}=0$ is an eigenvalue if and only if $\sum\,\chi_{i}(1)<|V|$ where the sum is over all characters with $\langle \tau,\chi_{i}\rangle\neq 0.$

Note, if $\Gamma=Cay(G,H)$ and $H$ is a normal connecting set of $G,$ then by \autoref{noor14} we have that each irreducible $G$-module say $U_i$ appears in exactly one eigenspace of $\Gamma$ with multiplicity of $m_i.$ Hence we conclude that $\langle \tau,\chi_i\rangle=\lambda_i m_i$ where $\lambda_i$ is an eigenvalue of $\Gamma$ with multiplicity of $m_i^2.$ \dne



Now, if $\Ga$ is singular, say $\gl_{1}=0,$ then every irreducible representation $\rho_{1,i}$ appearing in $\rho_{1}$ satisfies $\sum_{h\in H}\,\chi_{1,i}(h)=0,$ where $\chi_{1,i}$ is the character of $\rho_{1,i}.$ Conversely, if $\chi_{j,i}$ is an irreducible character with $\sum_{h\in H}\,\chi_{j,i}(h)=0$ then $\rho_{j,i}$ appears in $\rho$ and so there is some $E_{j}$ on which $\gl_{j}=0.$ \hfill $\Box$


\end{document}